\documentclass[12pt]{article}
\usepackage{latexsym}
\usepackage{amsmath, times, amsfonts, mathrsfs, amssymb}
\usepackage{amsmath}
\usepackage{amssymb}
\usepackage{graphicx}
\usepackage{subfigure}
\usepackage{natbib}
\makeatletter \oddsidemargin -0.1in \evensidemargin -.05in \textwidth
16cm \topmargin -1.8cm \textheight 25.3cm
\newcommand{\singlespacing}{\let\CS=\@currsize
\renewcommand{\baselinestretch}{1}\tiny\CS}
\newcommand{\doublespacing}{\let\CS=\@currsize
\renewcommand{\baselinestretch}{1.25}\tiny\CS}
\newtheorem{Theorem}{Theorem}[section]

\newcommand{\bea}{\begin{eqnarray}}
\newcommand{\eea}{\end{eqnarray}}
\newcommand{\be}{\begin{equation}}
\newcommand{\ee}{\end{equation}}

\newcommand{\bee}{\begin{eqnarray*}}
\newcommand{\eee}{\end{eqnarray*}}
\date{}
\begin{document}
\title{An inventory model with shortages for imperfect items using substitution of two products}
\date{}
\setcounter{page}{1}
\maketitle
\begin{abstract}
Inventory models with imperfect quality items are studied by researchers in past two decades. Till now none of them have considered the effect of substitutions to cope up with shortage and avoid lost sales. This paper presents an EOQ approach for inventory system with shortages and two types of products with imperfect quality by one way substitution. Our model provides significant advantage for substitution case while maintaining its inherent simplicity. We have provided numerical example and sensitivity analysis to justify the effectiveness of our model. It is observed that, presence of imperfect items affect the lot size of minor and major products differently. Under certain conditions, our model generalizes the previous existing models in this direction.
\end{abstract}
\hspace{1cm}\textbf{keywords-}{ inventory, EOQ, imperfect quality, substitution, shortages}

\section{Introduction}
\label{intro} Inventories are essentially required everywhere . From production to distribution and then to final customer, the role of inventory is incomparable. But absence of inventory irate customers and becomes catastrophic in the sense of future profit or goodwill.  The behavior of a customer is unpredictable when a product  becomes out of stock. The customer may often quit the store and get the same product from other places; or he may wait for the item; or he may take a similar product from the store itself. The phenomenon in the third case is known as product substitution.  For example, consider a customer who want to purchase a tourist bag of a particular brand. If it is not available in the shop the customer may be willing to purchase the other brand rather than going to other shop. The reason may be either shopkeeper insistence by offering similar product, or reputation of the shop or the other brand or it may be the customer's likeliness. This example demonstrates that how the inventory of a product will affect not only its own demand but also the demand of other products. Similar situations may frequently occur in the cases of food products, clothes, hard-drives, furniture,spare-parts, pharmaceuticals(drugs of various brands with identical composition). Therefore, a decision maker should consider the effect of substitution between products in determining the order quantities. In product substitution scenario, a collection of multiple products are available with the possibility to use an alternate product from the collection as a substitute whenever the original product runs out of stock.  In brief, substitution means demand for a certain quantity of a product is fulfilled by another product. Generally, there is a preferred product for satisfying a specific demand, which can optionally be substituted by certain alternative products (substitutes). The products between which substitution take place may be either different products manufactured by the organization itself or products replenished from one or more suppliers. Classical deterministic inventory models usually do not include product substitutions. However, neglecting substitution option in inventory systems may result in reduced efficiency in terms of customer satisfaction and costs. Sometimes in an inventory model, the products can be indexed in such a manner that a lower-index product may be used as substitute in place of a higher-index product. The product used to meet the demand of another product may be physically brought from another place and incur a transformation cost or it may be used in its original form. There are several advantages of substitution between products in inventory systems. Firstly, stockouts can by managed by the alternative manner using substitutes. It needs to store less amount of inventory of different variety of same type of product, leading to less amount of total holding cost.  By ordering larger lot sizes of a smaller number of products (that can substitute for other products), the total ordering cost and lead time can be reduced. Cheaper substitute is profitable in terms of purchase cost also. Substitutions can also be used to reduce the amount of outdated inventory of perishable products, e.g., by consuming substitute stocks first if they have an earlier expiry date. Also it may initiate revenue sharing contract among various systems which allows one to use other's inventory. This increases the possibility of inventory pooling to hedge against demand uncertainties and to help in reducing the safety stocks.  Sometimes this concept is also used as backup strategy, when the common product acts as a backup for the standard product. Demand may initially be satisfied by the standard product, and the common product will only be used to fulfill the demand if the standard product becomes out of stock.  Many cases are there when two or more generations (batch number, lot number or product version) of same product can be observed with different willingness of purchase by customers. This can also be seen as the two different products with different demands (e.g. Operating Systems like Windows7~and Windows8; Alpha and Beta versions of a software like MATLAB).  In this way, we can observe that substitution between products is very much practical and it is thus useful to include its effects on classical inventory models too.  Now we provide some literature survey which is related to our paper to grasp the matter more effectively.
\section{Literature Survey}
Classical Economic order quantity type models are used to determine the optimal  inventory policy when the demand is deterministic and the optimal ordering or production quantity are influenced by two opposite types of costs.  These models are used  to study the optimal inventory lot size that minimizes  inventory storage costs and maximizes systems' efficiency.  Although, these models have been successfully applied in the area of inventory management from early decades of previous century; yet they bear few unrealistic assumptions.  The fundamental EOQ model developed by \cite{Harris} involved the assumption of perfect quality items and no shortages ; both of these conditions fail to cope up with the realistic situations in the business scenario. In reality, the production process not always produce perfect quality items. Imperfect quality items are ineluctable in an inventory system due to imperfect production process, natural disasters, damages, or many other reasons. During last three decades, lot of research work were published in the area of EOQ and EPQ of imperfect quality items.  \cite{Rosenblatt} and \cite{Porteus} discussed deterministic inventory models where they assumed that the defective items could be reworked at a cost and found that the presence of a fraction of defective items motivate smaller lot sizes. \cite{Schwaller} presented a procedure and assumed that imperfect quality items are present in the lot in known proportions. He considered fixed and variable inspection costs for finding and removing the defective items. \cite{Zhang}considered a joint lot sizing and inspection policy with the assumption that a random proportion of lot size are defective. They assumed that defective items are not reworkable and used the concept of replacement of defective items by good quality items. \cite{Salameh} assumed that the defective items could be sold at a discounted price in a single batch by the end of the 100 \% screening process and found that the economic lot size quantity tends to increase as the average percentage of imperfect quality items increase. \cite{Goyal2} made some modifications in the model of \cite{Salameh} to calculate actual cost and actual order quantity. From a different viewpoint, \cite{Sheu} developed optimal production and preventive maintenance policy for an imperfect production system with  two types of out of control states.  \cite{Papachristos} pointed out that the sufficient conditions to prevent shortages given in \cite{Salameh} may not really prevent their occurrence and considering the withdraw time of the imperfect quality items from stock, they clarified a point not clearly stated in \cite{Salameh}. \cite{Wee1} developed the optimal inventory policy for an inventory system items with imperfect quantity where shortages were backordered. They allowed 100\% screening of items in which screening rate is greater than the demand rate. \cite{Chung} considered an inventory model with imperfect quality items under the condition of two warehouses for storing items. \cite{Jaber3} made certain assumption on learning curve and shown that percentage of defective lot size reduces according to the learning curve.  A detailed survey of the recent inventory models with imperfect items are provided by \cite{Khan}.  \cite{Roy} developed an EOQ model for imperfect items where a portion of demand were partially backlogged.  Amalgamating the effect of product deterioration, imperfect quality, permissible delay and inflation, \cite{Jaggi} investigated an inventory model when demand is a function of selling price.  \cite{Yu} extended  \cite{Salameh} when a portion of the defectives can be utilized as perfect quality and the utilization of the acceptable defective part will affect the consumption of the remaining perfect quality items in the stock. They have considered two cases for presence of  imperfect items as constant fraction and random fraction for their study.  In another model,  \cite{Jaber13} investigated Economic Order Quantity (EOQ) model for imperfect quality items  under the push-and-pull effect of purchase and repair option, when the option when defectives are repaired at some cost or it is replaced by good items at some higher cost. They presented two mathematical models; one for each case and discussed optimal policies for each case.

The effect of substitution in inventory systems has been studied by various researchers. In particular, the behavior of two-product substitution, where one product can substitute for the other, has been  investigated by numerous researchers pioneered by \cite{McGillivray}.  They assumed identical cost structures and substitution is stockout based but probabilistic. They obtained optimal policy using the simulation and heuristic approach. \cite{Parlar} considered two item inventory system with substitution under newsvendor framework  and obtained optimal order size. \cite{Drezner} were the first known researchers who considered the effect of shortage based substitution in deterministic inventory system for two types of items. They have shown that substitution is a profitable policy subject to a condition involving costs and demands, otherwise no substitution is optimal. They obtained the analytical expression of the optimal order quantities of both type items. \cite{Goyal1} suggested that it may be useful to consider time as a decision variable for the model of \cite{Drezner} and suggested an algorithm to obtain optimal policy. \cite{Gurnani} generalized the model of \cite{Drezner}  under hierarchical substitution scenarios. Since their formulation became intractable and was not providing the optimal solution, they have remodified their cost function in terms of optimal run out times and obtained the optimal inventory policy numerically. After that, no research work is available in the literature in deterministic inventory system with substitution. But many authors like  \cite{Nagarajan},\cite{Yadavalli},  \cite{Karakul} etc. have considered substitution in inventory control policy under newsvendor approaches. \cite{Pineyro} applied tabu search procedure to obtain maximizing remanufacturing quantity for economic lot-sizing problem with product returns and one-way substitution. \cite{Vaagen} considered  stochastic programming and simulation-based efficient optimization heuristics which can be applied to various models. They found that customer driven substitution are more costlier than retailer driven substitution. \cite{Liu} studied effect of substitution when two loss averse retailers are competing for substitutable product under stochastic demand rate and deterministic substitution rate. They have developed game theoretic model obtained unique optimal policy for Nash equilibrium under certain conditions. \cite{Tan} obtained the optimal service level for the profit maximizing policy of a retail product with Poisson arrival processes, stock-out based dynamic demand substitution, and lost sales using genetic algorithm. Recently, \cite{Salameh2} have solved joint economic lot size model with substitution and observed that substitution is highly beneficial in saving cost. Thus we observe that, previous researchers did not consider the effect of substitution in the EOQ with imperfect items in inventory management studies. The effect of substitute present in a deterministic inventory model can be used to reduce the inventory level. In the larger lot sizes due to imperfect production  a substitution effect comes into play that reduces the lot size which is the interesting part of this paper.  The main purpose of this present study is to develop a mathematical model with shortages for computing the economic order quantities where imperfect production process and substitution effect  are taken  into account.  \\
\\
~~~The organization of the paper is as follows. In Section 3, the basic model of \cite{Drezner} is reviewed and solved using time parameters.  In Section 4, we modified above model to formulate a new model that considers the effect of imperfect quality items and substitution. In Section 5, numerical example is presented and sensitivity analysis is performed to observe the effect of change of various parameters in the optimal policy. Section 6 gives conclusions and discusses various possible extension of this model. An appendix is also provided for some clarification.
\section{Brief review of the Drezner et al's model: Basic Model}
\cite{Drezner}  modified the EOQ model for two items by assuming a fraction of shortage of minor category item is met from major category item. Upon replenishing the inventory, both items are subject to face their own demands during the selling season. Minor items exhaust first and major items are delivered when the demand of minor items arises. The cycle ends when inventory of major items reaches to zero.  All other assumptions of this model remain to be the same as for the EOQ model; which include constant demand rates, no shortages, zero-lead time, etc. The behavior of inventory of minor and major items w.r.t. time  in \cite{Drezner} model is depicted in Figure - 1.\\
\vspace{0.1cm}
\hspace{1cm}
 In this model, it is assumed that first item is the major item and second item is the minor item which is substituted by the first item.  According to the requirement, if necessary, the first product is converted to second one, item-wise. ~Thus there is no holding cost for the second item after it is converted into first item.~There are two inventory cycles for each product shown in figure 1. We follow the hints provided by \cite{Goyal1} modified formulation of the model but explain it along the views of original model to adhere some level of simplicity.
\begin{center}
\includegraphics[width=15cm, height=10cm ]{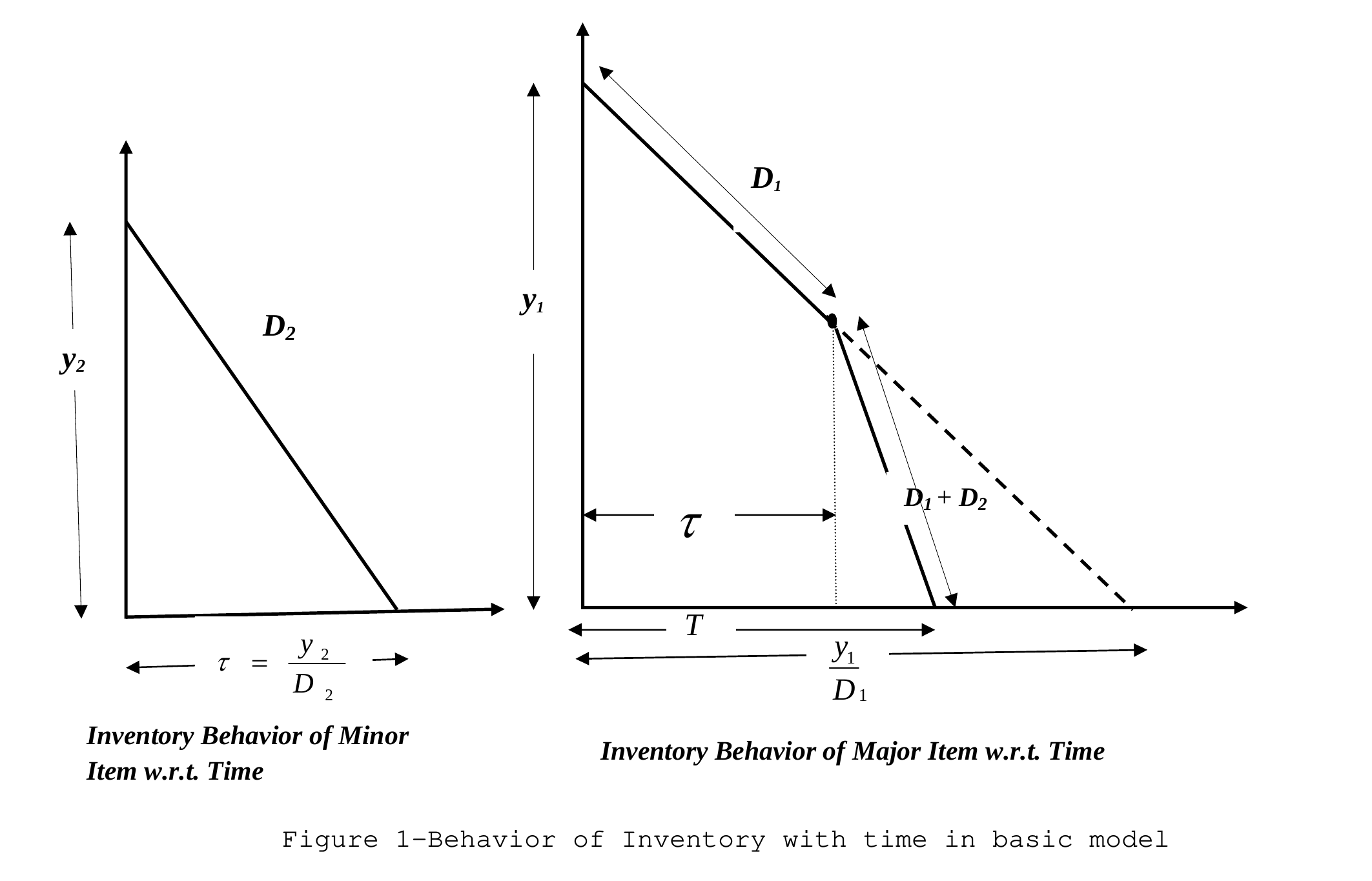}
\end{center}
Throughout this paper, we shall follow the notations given below. Let us assume following notations for reviewing the modified Drezner et al.'s\cite{Drezner} model:\\
$D_i=$ the annual demand for the product $i=1,2$;\\
$y_i=$ the order quantity for the product $i=1,2$;\\
$c_{hi}=$ the holding cost for the product $i=1,2$;\\
$c_t=$ the transfer cost for product into product 2 \\
$c_o=$ the ordering cost for both product \\
$\tau=$ the time interval during which no substitution is required.\\
$T=$ cycle time.
\par Now, holding cost of each type items:\\
$H_1=c_{h1}\big[\int_0^{T} (D_1 t) dt+ \int_0^{\tau}(D_2 t) dt\big]$\\
\vspace{0.3cm}
\hspace{0cm}
$=c_{h1}\big[D_1\frac{ T^2}{2}+D_2\frac{(T^2-{\tau}^2)}{2}\big]$
\\
and,
$H_2=c_{h2}\big[\int_0^{\tau} D_2 t dt\big]=c_{h2} D_2 \frac{{\tau}^2}{2} $ \\  \vspace{0.3cm}

Transfer cost of substituted item: \\  \vspace{0.3cm}
$CT=D_2 c_t (T-\tau)$  \\  \vspace{0.3cm}
Total average cost of system:
\begin{eqnarray*}
TAC=\frac{c_o+H_1+H_2+CT}{T}
\end{eqnarray*}
\begin{eqnarray}
TAC=\frac{c_{o}}{T}+c_{h1}[\frac{D_1}{2}T + \frac{D_2}{2}(T-\frac{\tau^2}{T})]+c_{h2}\frac{D_2}{2} \frac{\tau^2}{T}+D_2  c_t (1-\frac{\tau}{T})
\end{eqnarray}
In above formulation there will be three cases :\\ \vspace{0.1cm} \hspace{1.5cm}
(i) \emph{Partial Substitution},~~~ (ii)  \emph{Full Substitution}, ~~and~~ (iii) \emph{No Substitution}. \\ \vspace{0.1cm} \hspace{.5cm}
We first consider the case of \textbf{partial substitution case} in which $0<\tau<T$ and $y_2>0$.\\
We have from $\frac{\partial TAC}{\partial \tau}=0$ ,
  \begin{eqnarray}
\frac{-c_{h1} D_2 \tau+c_{h2} D_2 \tau - D_2 c_t }{T}=0
  \end{eqnarray}
which gives the optimal value of $\tau$ as,
\begin{eqnarray}
{\tau}^*=\frac{c_t }{c_{h2}-c_{h1}}
\end{eqnarray}
Also,~  $\frac{{\partial}^2 TAC}{\partial {\tau}^2}=D_2 (c_{h2}-c_{h1}) $ \\
\vspace{0.5cm}
which is positive since it is assumed that $c_{h2}>c_{h1}$. Thus $TAC$ is convex. \\
\vspace{0.2cm}
But if $\frac{c_t }{c_{h2}-c_{h1}}>T$, we observe that time interval during which no substitution is necessary, is more than the cycle time $T$ which is absurd. Thus we have $0\leq \tau \leq  T$ for all practical cases.\\
 \\
 ~~~Now from $\frac{\partial TAC}{\partial T}=0$~, we obtain:
 \begin{eqnarray}
  c_{h1} \bigg(\frac{D_1+D_2}{2}+\frac{D_2 {\tau}^2}{2 T^2}\bigg)-c_{h2}\frac{ D_2 {\tau}^2}{2 T^2}+D_2 c_t\frac{ {\tau}}{T^2}-\frac{c_o}{T^2}=0
\end{eqnarray}
From equation (4), the optimal value of $T^*$ is obtained as:
 \begin{eqnarray*}
  T^*=\sqrt{\frac{c_o+\frac{D_2 \tau^2}{2} (c_{h2}-c_{h1})-D_2 {c_t}}{c_{h1}(\frac{D_1+D_2}{2})}}
\end{eqnarray*}
Substituting the value of $\tau$ from equation (3), we obtain $T^*$  as follows:
  \begin{eqnarray}
  T^*=\sqrt{\frac{2 c_o-\frac{D_2 {c_t}^2}{c_{h2}-c_{h1}}}{c_{h1}(D_1+D_2)}}
\end{eqnarray}
and the corresponding  optimal order quantity $y_2^*$ becomes:\\   \vspace{0.2cm}
   \begin{eqnarray}
   y_2^*={\tau}^* D_2=\frac{c_t D_2 }{c_{h2}-c_{h1}}
\end{eqnarray}
The optimal order quantity of type-1 items  $y_1^*$ is clearly the required order quantity for satisfying demand for both items during whole season $T_1$ minus the optimal order quantity  ${y_2^*}$.
      \begin{eqnarray}
   y_1^*={T}^*(D_1+ D_2)-y_2^*={T}^*(D_1+ D_2)-\frac{c_t D_2 }{c_{h2}-c_{h1}}
\end{eqnarray}
where $T^*$ is given by (5).
{\Theorem There exist unique optimal substitutional economic order quantity for all possible parameter values if $c_t<\sqrt{\frac{2c_o(c_{h_2}-c_{h_1})}{D_2}}$.}\\
\noindent \textbf{ Proof : } We know that for a convex function local optimum is global optimum. In order to prove the unique optimal solution, it suffice to show the convexity of total cost function.
Let us re-write the equation (1) as, \\
 $ TAC=\frac{c_{o}}{T}+c_{h1}\big[\frac{D_1}{2}T + \frac{D_2}{2}(T-\frac{\tau^2}{T})\big]+c_{h2}\frac{D_2}{2} \frac{\tau^2}{T}+D_2  c_t (1-\frac{\tau}{T})$
 \\
\vspace{0.1cm}
Simple calculation gives:
\\
\vspace{0.1cm}
 $\frac{\partial^2 TAC(\tau, T)}{\partial\tau^2}=\frac{D_2(c_{h2}-c_{h1})}{T}$\\
\vspace{0.1cm}
$\frac{\partial^2 TAC(\tau, T)}{\partial T^2}=\frac{1}{4T^3}\big[2 c_o-\frac{D_2 c_t^2}{c_{h2}-c_{h1}}\big]$\\
\vspace{0.1cm}
$\frac{\partial^2 TAC(\tau, T)}{\partial\tau \partial T }=D_2\frac{c_t-(c_{h2}-c_{h1})}{T^2}$
\\
\vspace{0.1cm}
Now we have by definition,\\
\vspace{0.2cm}
\hspace{1cm}
Hessian $H_0(\tau,T)=\frac{\partial^2 TAC(\tau, T)}{\partial\tau^2}\frac{\partial^2 TAC(\tau, T)}{\partial T^2}-\Big(\frac{\partial^2 TAC(\tau, T)}{\partial\tau \partial T }\Big)^2$\\
$=\frac {D_2 \left( \frac{1}{2}\,(c_{h2}-c_{h1})c_o- \left(  \frac{3}{4}\,{(c_{h2}-c_{h1})}^{2}{\tau}^{2}- \frac{3}{2}\,(c_{h2}-c_{h1})\tau\,c_t+{c_t}^{2} \right) D_{{2}} \right) }{T^4}$  \\
Putting the value of $\tau$ from (3), we obtain:\\

$H_0(\tau,T)=\frac{2D_2}{T^4}\big[c_o(c_{h2}-c_{h1})-\frac{1}{2}D_2 c_t^2\big]$ \\
\vspace{0.3cm}
\hspace{0.5cm}
Clearly, $H_0(\tau,T)>0$ if
\vspace{0.1cm}
\hspace{0.5cm}  $c_t<\sqrt{\frac{2c_o(c_{h_2}-c_{h_1})}{D_2}}$ \\
\vspace{0.1cm}
\hspace{0.5cm}
Hence the proof follows.\\
\\
Now we study the special cases of above model to incorporate the situations when all the items of minor products are substituted completely by major products and also the case when substitution is not required.
\begin{center}
\textbf{ SPECIAL CASES}\\
\end{center}
\textbf{Case 1: Full Substitution} :\\
\vspace{0.3cm}
In this case, $\tau=0$,  $y_2=0$ and $T=\frac{y_1}{D_1+D_2}$
 \\
\vspace{0.1cm}
We have, holding cost of major items $H_1=\frac{c_{h1}}{2}(D_1+D_2)T^2$, and\\
holding cost of minor items $H_2=0$
 \\
\vspace{0.1cm}
Transformation cost $CT= c_t (D_2 T)$.
 \\
\vspace{0.1cm}
Now the total cost becomes,
  \begin{eqnarray}
  TAC=\frac{c_{h1}}{2}(D_1+D_2)T +D_2 c_t+\frac{c_o}{T}
    \end{eqnarray}
 $\frac{\partial TAC}{\partial T}=0$ gives,
  \begin{eqnarray}
T^*=\sqrt{\frac{2 c_o}{c_{h1}(D_1+D_2)}}
  \end{eqnarray}
  \\
\vspace{0.1cm}
Also  $\frac{\partial^2 TAC}{\partial T^2}=\frac{2 c_o}{T^3}>0$ \\
   \\
\vspace{0.3cm}
Thus we have the only optimal order quantity of major item:
\begin{eqnarray}
y_2=D_2 T=D_2 \sqrt{\frac{2 c_o }{c_{h1}(D_1+D_2)}}
\end{eqnarray}
\textbf{Case 2: No Substitution} :\\
\vspace{0.1cm} \hspace{0.01cm}
In this case we have, $\tau=T$~ and $\frac{y_1}{D_1}=\frac{y_2}{D_2}=T$\\
\vspace{0.1cm}
The holding costs of major items and minor items become \\
\vspace{0.1cm} \hspace{0.01cm}
$H_1=\frac{c_{h1} D_1 T^2}{2}$,~~$H_2=\frac{c_{h2} D_2 T^2}{2}$~~~respectively.\\
\vspace{0.1cm}
Also, transformation cost~$CT=0$\\
\vspace{0.1cm}
Now the total cost becomes,
  \begin{eqnarray}
TAC=(c_{h1} D_1+c_{h2} D_2)\frac{T}{2}+\frac{c_o}{T}
  \end{eqnarray}
First order condition~ $\frac{\partial TAC}{\partial T}=0$ gives,
   \begin{eqnarray}
 T^*=\sqrt{\frac{2 c_o}{c_{h1}D_1+c_{h1} D_2}}
 \end{eqnarray}
In this case, also  $\frac{\partial^2 TAC}{\partial T^2}=\frac{2 c_o}{T^3}>0$~\\
-which means that convexity criteria of total average cost$TAC$ is satisfied.\\
\vspace{0.2cm}
Thus we have:\\
 the optimal order quantity of minor item:\\
\vspace{0.1cm} \hspace{5cm} $y_1=D_1 T^*=D_1 \sqrt{\frac{2 c_o}{c_{h1}D_1+c_{h1} D_2}}$  \\
 the optimal order quantity of major item:\\
\vspace{0.1cm} \hspace{5cm} $y_2=\frac{D_2 }{D_1}y_1=D_2 \sqrt{\frac{2 c_o}{c_{h1}D_1+c_{h1} D_2}}$\\
\vspace{0.4cm}
It can be observed that above formulae are identical with corresponding results in \cite{Drezner} model.
\section{Our Mathematical Model: EOQ model for imperfect items with shortage and substitution(EOQISS)}
\begin{center}
\includegraphics[width=14cm, height=15cm ]{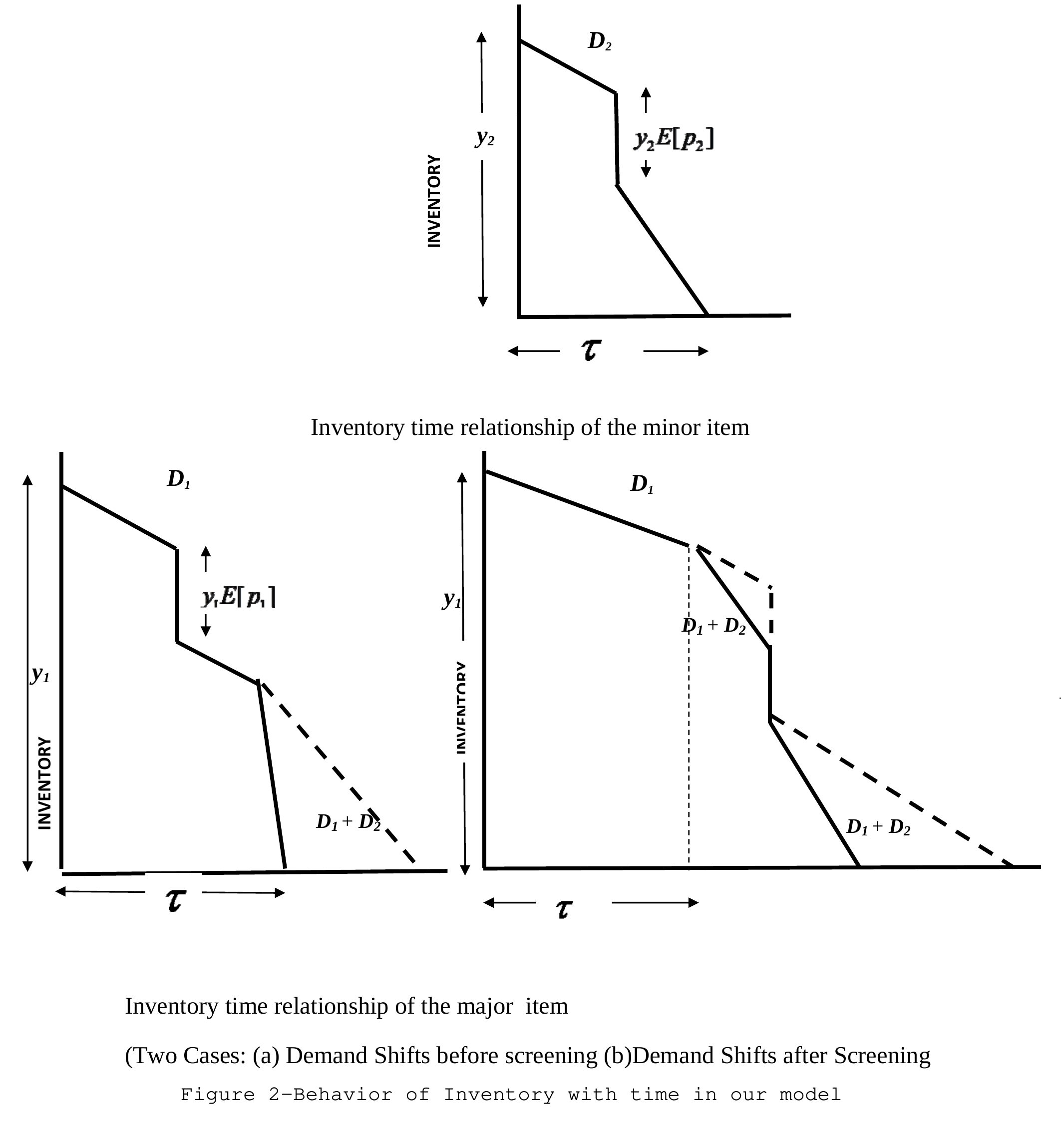}
\end{center}
In this section, we incorporate the presence of imperfect items  to the problem discussed in aforementioned model.  The behavior of inventory of minor and major items w.r.t. time  in our model is shown in Figure - 2. We need few more notations and assumptions to develop our model which are explicitly provided in forthcoming subsections.
\subsection{Notations}
\label{theory} In order to incur such modification we use following additional notations :\\
$p_i=$ the fraction of imperfect items in the lot of the product $i=1,2$;\\
$t_{s i}=$ the time epoch at which the screening finishes for the product $i=1,2$;\\
$x_{ i}=$ the screening (rate)speed for the product $i=1,2$;\\
$E[\cdot]$ denote probabilistic expectation operator.\\
We have $t_{s i}=\frac{y_i}{x_i}$

\subsection{Assumptions}
The following assumptions are used to establish our mathematical model:
\begin{enumerate}
\item Replenishment is instantaneous, Lead time is negligible for both type of items.
\item The demand rate of each type of items are constant and unequal in general.
\item Substitution of first item by second item is allowed only when first item is out of stock (Stockout based substitution).
\item Substitute is abundant whose shortage never occur.
\item A fraction of both types of items received contain imperfect items; which are screened before sending to customers. Such fraction are random variables whose probability distribution follow uniform  distribution.
\item Screening rates of both type of defective items greater than their demand rates $x_i>D_i$ for $i=1,2$
\item The following relation holds between expected value of defective items:\\
~~~ $E[p_i]<1-\frac{D_i}{x_i}$; ~~where $i=1,2$.
\item  The inventory reaches zero level at the end of each cycle.
\item The following relation holds between holding cost parameters:\\
$c_{h2}>c_{h1}$
\end{enumerate}
The model we are developing are subjected to two scenarios:\\
\vspace{.1cm} \hspace{5.3cm}
(1). $ \tau \le t_{s 1}$  ~~~~~and~~~~~~
(2). $ \tau \geq t_{s 1}$ \\
\vspace{.1cm} \hspace{.5cm}
Note that in both cases, number of perfect quality items of major items remain same. This imply when shortage of minor items will occur, it will be meet from perfect quality major type items.  Thus in effect, there is no change in the corresponding optimal policies. \\
\vspace{.2cm} \hspace{0.1cm}
The differential Equations governing the inventory level at time $t\in [0,t_{s 2}]$
\begin{eqnarray}
\frac{dI_2}{dt}=-D_2
\end{eqnarray}
subject to the initial condition:~~$I_2(0)=y_2$.\\
\hspace{.1cm} Solving and applying the initial condition,\\
$I_2(t)=y_2-D_2 t$ \\
At time  \hspace{.05cm} $t_{s 2}$,  the initial inventory~ $=y_2-D_2 t_{s 2}$\\
At time  \hspace{.05cm} $t_{s 2}$, the final inventory ~$I_{s2}=y_2(1-E[p_2])-D_2 t_{s 2}$\\
The differential equation governing the inventory level at time $t\in [t_{s 2},\tau]$ is\\
\begin{eqnarray}
\frac{dI_2}{dt}=-D_2
\end{eqnarray}
subject to $I_2(t_s)=I_{s2}$  \\
\vspace{.1cm}
The differential equation governing the inventory level at time $t\in [0,\tau]$ is
\begin{eqnarray}
\frac{dI_1}{dt}=-D_1
\end{eqnarray}
subject to the initial condition:~~$I_1(0)=y_1$.\\
\vspace{.1cm}
\hspace{.1cm}
Solving and applying the initial condition, we obtain:
\begin{eqnarray}
I_1(t)=y_1-D_1 t
\end{eqnarray}
The differential equation  governing the inventory level at time $t\in [\tau, t_{s1}]$ is
\begin{eqnarray}
\frac{dI_1}{dt}=-(D_1+D_2) \forall t\in [\tau, t_{s1}]
\end{eqnarray}
Solving (17), we get\\
\vspace{0.1cm}
\hspace{1cm}$I_1(t)=y-D_1 t -D_2( t-\tau)$\\
\vspace{.2cm}
At time \hspace{.05cm} $t_s$, the final inventory \hspace{.05cm} $I_{s1}=y_1(1-E[p_1])-D_1t_{s1}-D_2( t_{s1}-\tau)$\\
\vspace{.2cm}
The differential equation governing the inventory level at time $t\in (t_s,T]$
\begin{eqnarray}
\frac{dI_1}{dt}=-(D_1+D_2) \forall t\in [t_{s1},T]
\end{eqnarray}
\vspace{.1cm}
\hspace{.2cm} subject to the initial condition:\hspace{1.15cm}$I_2(t_s)=I_{s1}$.\\
\vspace{.1cm}
Solution of (18) is\\
\vspace{.1cm}
\hspace{2cm}
 \begin{eqnarray}
 I_1(t)=y_1(1-E[p_1])-D_1 t+D_2 (t-\tau)
 \end{eqnarray}

Now the various cost items are calculated as follows:\\
\vspace{0.1cm}
\hspace{0.5cm}
$OC=c_o$\\
\vspace{0.1cm}
\hspace{0.5cm}
$HC=c_{h2}\int_0^{\tau}I_2(t) dt+c_{h1}\int_0^TI_1(t) dt$~~~and~~~$CT=c_t D_2 (T-\tau)$\\
\vspace{0.1cm}
\hspace{0.5cm}
The expression of total cost is calculated as:
\begin{eqnarray}\label{mm1}
\nonumber    TC =&& c_{h1} \Big[ \frac{ D_ 1 +D_ 2 }{2} {T}^{2}+\frac { ( ( T-\tau  ) D_ 2 +T D_ 1 )^2 E[p_1] }{  (1- E[p_ 1] ) ^{2} x_ 1 }  -\frac{1}{2} D_ 2 {\tau}^{2} \Big] +
\\ && c_ { h2} \Big[\frac{1}{2}D_ 2 {\tau}^{2}+\frac{E[p_ 2] D_ 2 {\tau}^2}{(1- E[p_ 2] )^{2} x_ 2}\Big] +c_ o +D_ 2 c_ t ( T-\tau  )
\end{eqnarray}
We have, Total Average Cost~ $TAC= \frac{TC}{T}$\\
\vspace{0.1cm}
Hence,
\begin{eqnarray}\label{mm99}
\nonumber   TAC&=&c_{h1}\Big[\frac{D_1+D_2}{2} T+{\frac { ( ( T-\tau ) D_2+T D_1 ) ^{2}E[p_1]}{ ( 1-E[p_1] ) ^{2}x_1T}}-{\frac {D_2{\tau}^{2}}{2 T}}
 \Big] +c_{h2} \Big[ {\frac {D_2{\tau}^{2}}{2 T}}+{\frac {E[p_2] D_2{\tau}^{2}}{ ( 1-E[p_2] ) ^{2}x_2T}} \Big] \\ && +{\frac {c_o}{T}}+D_2 c_t \Big[1-{\frac {\tau}{T}}\Big]
\end{eqnarray}
 \vspace{.1cm}
The convexity of $TAC$ guarantees the unique optimal solution in terms of $\tau$ and $T$. \\
We now state the following theorem for convexity of cost function as well optimality and uniqueness of the optimal solution:
{\Theorem If the parameters of the cost function given by equation (21) satisfy \\
 $c_o>\frac{c_t^2}{2(c_{h2}-c_{h1})}+\tau(T-1)c_t D_2$ , then the optimal cost can be obtained uniquely with optimally unique parameters $\tau^*$ and $T^*$.}
\\
\vspace{0.4cm} For Proof, see \textbf{Appendix-I}.\\
\vspace{0.4cm}
 We obtain first order condition for $\tau$, by differentiating $TAC$ w.r.t. $\tau$, and equating to zero, that is,~  $\frac{\partial TAC}{\partial \tau}=0$, which gives,
\begin{eqnarray}
\nonumber c_{h1} (-2\,{\frac { ( (T-\tau ) D_2+TD_1 ) E[p_1]D_2}{ (1-E[p_1] ) ^{2} x_1 }}-D_2\tau ) +c_{h2} (D_2\tau+2\,{\frac {E[p_2] D_2\tau}{ (-(1-E[p_2]) ) ^{2}x_2}} ) -D_2 c_t =0 \\
\end{eqnarray}
Simplification of equation(22) gives
\begin{eqnarray}
\nonumber \tau=2\,{\frac { ( 1-E[p_2] ) ^{2} ( \frac{1}{2} c_{t}  ( 1-E[p_1] ) ^{2} x_1+T c_{h1} E[p_1]( D_{{1}}+D_{{2}}  )  )x_2}{ (  ( 1-E[p_2] ) ^{2} ( c_{h2}-c_{h1} )x_2+2\,c_{h2}\, E[p_2]  )  ( 1-E[p_1] ) ^{2}x_1+2\, c_{h1} x_{{2}} D_{{2}} E[p_1]( 1-E[p_2] ) ^{2}}} \\
\end{eqnarray}
Substituting $E[p_1]=E[p_2]=0$ in equation (23), we obtain,\\
\vspace{0.1cm}
\hspace{5.2cm} $ \tau=\frac{c_{t} }{c_{h2} -c_{h1} }$ \\
\vspace{0.1cm}
which is same as equation(3).
\\
\vspace{0.1cm}
Hence the optimal order quantity of minor items $y_2^*$ is given by:
\begin{eqnarray}
y_2^*=\frac{D_2 \tau}{1-E[p_2]}
\end{eqnarray}
where $\tau$ is given by equation (23).
\\
In order to obtain optimal order quantity of major item $y_1$, we require to find the expression for the optimal value $T^*$. \\
We differentiate (21) w.r.t. $T$, and equating to zero, we obtain the value of $T^*$ as:
\begin{eqnarray}
T^*=\frac{\sqrt{A_1}}{c_{h1} x_2 ( D_1+D_2 )  ( 1-E[p_2]  )  ( \frac{x_1}{2} (1-E[p_1] )^{2}+ ( D_1+D_2 ) E[p_1] )}
 \end{eqnarray}
where $A_1$ is given by:\\
$A_1=c_{h1} x_2  ( D_1+D_2  )   ( -\frac{1}{2}  ( 1-E[p_1]  ) ^{2}  (   ( \tau  (   ( c_{h1}-h_{
{2}}  ) \tau+2 c_t  ) D_2- 2 c_o  )
  ( 1-E[p_2]  ) ^{2}x_2+2\,c_{h2}D_2{\tau}^{2} E[p_2]  ) x_1+c_{h1}x_2E[p_1]{D_2}^{2}{\tau}^{2}
  ( 1-E[p_2]  ) ^{2}  )   ( \frac{1}{2}  ( 1-E[p_1] ) ^{2} x_1+  ( D_1+D_2  ) E[p_1] )$ \\  \vspace{0.1cm}
\hspace{1cm} Now the optimal order quantity of major item is given by:
\begin{eqnarray}
y_1=\frac{( D_1+D_2 ) T^*-D_2\tau^*}{1-E[p_1]}
 \end{eqnarray}
 The equation (26) can be visualized as the extension of equation (5) in terms of imperfect quality items. \\
\vspace{.1cm}
\hspace{1cm} By substituting $E[p_1]=E[p_2]=0$ and $\tau=\frac{c_t}{c_{h2}-c_{h1}}$, we obtain $T^*$ which corresponds to the same value as in equation $(5)$.

 \begin{center}
\textbf{ SPECIAL CASES}\\
\end{center}
\textbf{Case 1: Full Substitution} :\\
\vspace{0.3cm}
In this case, $\tau=0$,  $y_2=0$ and $T=\frac{y_1(1-E(p_1))}{D_1+D_2}$
 \\
\vspace{0.1cm}
We have $H_1=\frac{c_{h1}}{2}(D_1+D_2)T^2$~,~ $H_2=0$~~~ and~~~$CT=D_2 c_t T$ .\\
\vspace{0.1cm}
Thus we obtain
\begin{eqnarray}
\nonumber TAC&=&c_{h1} ( \frac{1}{2}\, ( D_1) T+{\frac {T D_1 ^{2} E[p_1]}{ ( 1-E[p_1] ) ^{2}x_1}} ) +c_{h2} \big( \frac{1}{2}\,D_2{T}+{\frac {E[p_2] D_2{T}}{ (1-E[p_2] ) ^{2}x_2T}} \big) +{\frac {c_o}{T}}
\\&&
 \end{eqnarray}
 Above cost function consists of single variable $T$ only and it is obviously convex because: \\
\vspace{0.1cm}
  $\frac{\partial^2 TAC(\tau, T)}{\partial T^2}=\frac{2c_o}{T^3}>0$ \\
\vspace{0.1cm}
Now the first order optimality condition w.r.t. variable $T$ ,~that is,  $\frac{\partial TAC}{\partial T}=0$ gives,
  \begin{eqnarray}
T^*=\sqrt{\frac{2 c_o x_1 x_2 ( 1-E[p_1] )^{2} ( 1-E[p_2] )^{2}}{x_1 ( 1-E[p_1] )^{2}\big((c_{h1}D_1+c_{h2}D_2)( 1-E[p_2] )^{2}+2c_{h2}D_2 E[p_2]\big)+2 x_2 c_{h1}D_1^2 E[p_1]( 1-E[p_2] )^{2}}}
  \end{eqnarray}
Now, we obtain the value of optimal order quantity of major item $y_1^*$ using equation (26) by substituting the value of $T^*$ calculated above and putting $\tau^*=0$ in equation (26).
 \\
\vspace{0.2cm}
Thus~~ $y_1^*=\frac{(D_1+D_2)T^*}{1-E[p_1]}$
 \\
\vspace{0.3cm}
\textbf{Case 2: No Substitution} :\\
\vspace{0.1cm} \hspace{0.01cm}
In this case we have, $\tau=T$.~\\
Total average cost function becomes
\begin{eqnarray}
\nonumber TAC&=&c_{h1} ( \frac{1}{2}\, ( D_1+D_2 ) T+{\frac {T ( D_1+D_2)^{2}E[p_1]}{ (1-E[ p_1] ) ^{2} x_1}} )+{\frac {c_o}{T}}+D_2 c_t
 \end{eqnarray}
  In this case also, the convexity of $TAC$ holds because  $\frac{\partial^2 TAC}{\partial T^2}=\frac{2 c_o}{T^3}>0$  \\
\vspace{0.1cm}
The first order condition gives: $\frac{\partial TAC}{\partial T}=0$ gives,
   \begin{eqnarray}
 T^*=\sqrt{\frac{2 c_o x_1}{c_{h1}(D_1+ D_2)\big((D_1+D_2)E[p_1]+\frac{x_1(1-E[p_1])^2}{2}\big)}}
 \end{eqnarray}
Now, from equation (24) using the value of $T^*$ calculated above and $\tau^*=T^*$, we derive the value of $y_1^*$ as \\
\vspace{0.1cm}
 $y_1^*=\frac{D_1 T^*}{1-E[p_1]}$\\
\vspace{0.1cm}
Also we have from equation (24)\\
 $y_2^*=\frac{D_2 T^*}{1-E[p_2]}$\\
 \vspace{0.1cm}
 \hspace{1cm} Above expression of $y_1^*$ and $y_2^*$ shows that each of the major and minor items satisfies their demand independently and there is no demand interactions- which is the consequence of no substitution assumption.\\
\textbf{ Remark:}\\
 Since some moderate cost is associated with the transformation of minor item to major item and the holding cost of major item is higher than that of minor item; it becomes obvious that full substitution is never optimal. This was also true for \cite{Drezner} model.  This will help us to develop  simple solution procedure which is provided in next subsection.
\subsection{Solution Procedure}
  Due to suboptimality of full substitution, we have developed an algorithm resembling the algorithm of \cite{Goyal1} that will determine the optimal policy of our model for the case of partial substitution and no substitution. The algorithm is straightforward and provides optimal policy in simple way. We now state the algorithm as follows: \\
\\
\textbf{Algorithm:}\\
\\
\emph{\textbf{Step-1}}: Calculate $T_0=\sqrt{\frac{2 c_o x_1}{c_{h1}(D_1+ D_2)\big((D_1+D_2)E[p_1]+\frac{x_1(1-E[p_1])^2}{2}\big)}}$ ~~~and\\
 $\tau^*=2\,{\frac { ( 1-E[p_2] ) ^{2} ( \frac{1}{2} c_{t}  ( 1-E[p_1] ) ^{2} x_1+T c_{h1} E[p_1]( D_{{1}}+D_{{2}}  )  )x_2}{ (  ( 1-E[p_2] ) ^{2} ( c_{h2}-c_{h1} )x_2+2\,c_{h2}\, E[p_2]  )  ( 1-E[p_1] ) ^{2}x_1+2\, c_{h1} x_{{2}} D_{{2}} E[p_1]( 1-E[p_2] ) ^{2}}} $
 \\
\\
\emph{\textbf{Step-2}}:  If $T_0>\tau$, go to next step (partial substitution case); otherwise, assign the optimal values(no substitution case) as follows:\\
$T^*=T_0$ ~~~ $y_1^*=\frac{(D_1)T^*}{1-E[p_1]}$~~ $y_2^*=\frac{(D_2)T^*}{1-E[p_2]}$ and stop.
\\
\\
\emph{\textbf{Step-3}}: Evaluate optimal values of $T^*$ using equation (25).\\ With this value of $T^*$, optimal order quantities are then obtained as~
$y_1^*=\frac{(D_1+D_2)T^*-D_2\tau^*}{1-E[p_1]}$ ~~and~~$y_2^*=\frac{(D_2)\tau^*}{1-E[p_2]}$.\\
\section{Numerical Example and Sensitivity Analysis}
\subsection{Numerical Example}
In order to compare our work with that of \cite{Drezner}, we have taken the data from their paper and solved the basic model in terms of $\tau$ and $T$ and obtained identical result for order quantities.
The parameters are as follows:\\
\vspace{0.1cm} \hspace{0.2cm}
 $ D_1 = D_2 =1000$ units/year, ~ $c_{o}$ =  \$4500/cycle and $c_{h1}$ = \$1/unit/year.\\
 We include following new parameters in our model
 $x_1 = 175 200$, $x_2 = 175 100$,\\
 Let $p_1$ and $p_2$ be uniformly distributed random variables whose expected values are $E[p_1]$ and $E[p_2]$.
 We consider following three cases:
 (a) $E[p_2]-E[p_1]$ is very small, (b)$E[p_2]-E[p_1]$ is high. (c) $E[p_2]<E[p_1]$
 Beside that, we consider the cases of partial substitution, full substitution and no substitution situations in basic model as well as our (EOQISS) model.

 Numerical results of basic model which matches with \cite{Drezner} model is given in Table 1.
 Results for our(EOQISS) model for three different pair of values of $E[p_1]$ and $E[p_2]$  is shown in Table 2.

\begin{center}
\includegraphics[width=15.5cm, height=4.3cm ]{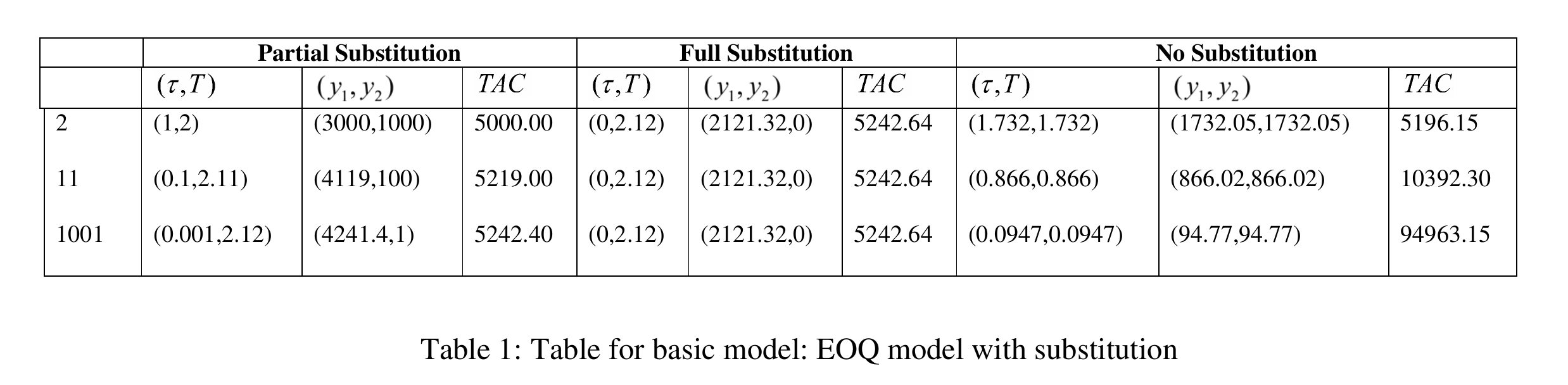}
\end{center}
 \begin{center}
\includegraphics[width=18cm, height=13.5cm ]{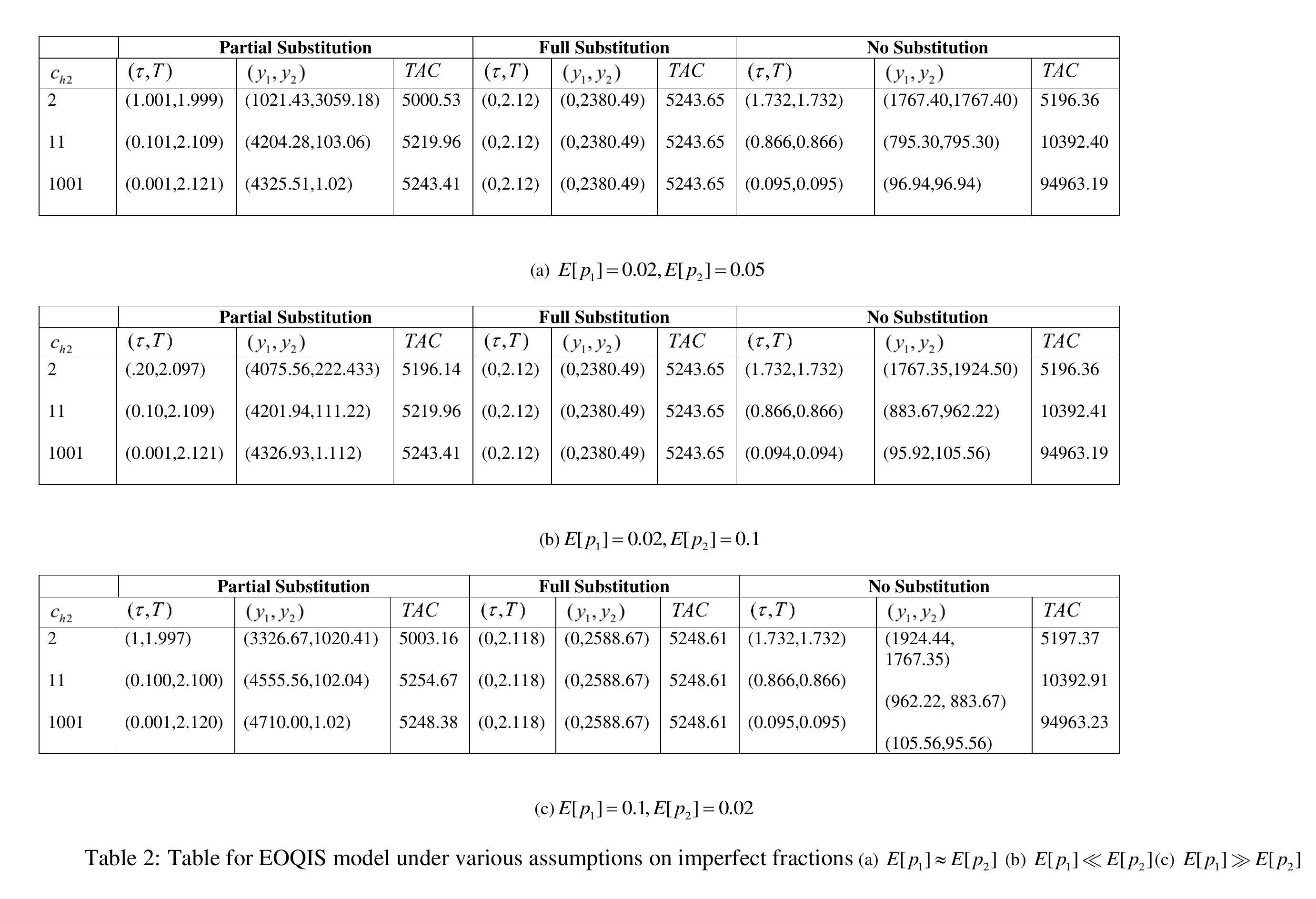}
\end{center}

From the above tables, we infer following observations for the inventory model with imperfect quality items:
\begin{itemize}
\item For the partial substitution case, we observe that:
 \begin{itemize}
 \item The order quantity $y_2$ of the minor item obtained by our model is very nearly equal to the corresponding order quantity in the basic model without imperfect items if the expected value of imperfect fraction of minor item is very less. But with small difference w.r.t. expected imperfect fraction value of the major item, the order quantity of major item increases moderately when holding cost of the second item increases gradually.
 The total average cost in above situation is very nearly equal to the corresponding total average cost for perfect items because we have assumed that screening cost is negligible.
\item Order quantity $y_2$ of the minor item obtained by our model decreases when its lot contains higher expected fraction of defective items. This is quite logical in the sense that it reduces the total expected number of defective items in the inventory. Consequently the order quantity of the major item increases and total cost also increases.
     \item Order quantity $y_1$ of the major item obtained by our model increases rapidly than that of minor items when its lot contains higher expected fraction of defectives.  The increase in the lot size of minor item is very less (about 2\% in our case). Corresponding costs have also minor variation compared to basic model.
    \end{itemize}
\item For the case of no substitution, we observe that time at which shortages of minor and major item occur are their corresponding cycle times. Both items are satisfied by their own demand and have independent inventory policies. Also we notice the following:
   \begin{itemize}
   \item  When the expected fraction of defectives in both type of items are small, the order quantity of each item increases as the holding cost of the second items increases gradually.~Corresponding optimal cost increases negligibly.
      \item  When the expected fraction of minor items' defectives are higher than that of major items, the optimal lot size increases moderately but the corresponding optimal cost increase negligibly.
       \item   When the expected fraction of major items' defectives are higher than that of minor items, the optimal lot size increases moderately but the cost increase marginally.
       \end{itemize}
\item For the case of full substitution, we observe that order quantity $y_2=0$ and run-out time $\tau=0$.  But the lot sizes and costs remain identical for variation in holding costs. Since there is no minor item, therefore $E[p_2]=0$. Also we notice that:
   \begin{itemize}
   \item  When the expected fraction of major items' defectives  increases, the order quantity of major items increases rapidly corresponding to other cases. Associated optimal cost also increases rapidly.
   \item   Thus we can say that presence of defective items play major role in full substitution case and order quantity and optimal cost are most sensitive to  expected fraction of defectives.
       \end{itemize}
    \end{itemize}
It can be stated from the above discussion that other than holding costs, the presence of defective items also play major role in the decision making regarding optimal inventory and costs in all cases of our model.

\section{Conclusion}
In this paper, we attempt to investigate a two-product deterministic inventory system with imperfect quality items under shortages and one way substitution. We have used an effective method which considers time period viz. run-out time and cycle time as decision variables. Using this approach, we have obtained a simplified closed form expression for the optimal order quantities which is comparable with the model of \cite{Drezner}. Our work differs from previous research works in the sense that we consider the presence of a fraction of imperfect items in the system, which is screened out and sold at some negligible salvage price at the end of the selling season. Thus, the contribution of this paper is two-fold. First, we modify the condition for optimality (in terms of cycle time) of cost function of the basic model of \cite{Drezner} without imperfect items. Another contribution of our paper is to incorporate the effect of imperfect quality items in both type of items when the substitution of a product is permitted. We observed the different behavior of the optimal order quantities and corresponding optimal cost when expected defective fraction is more in major item and otherwise. In particular, it is observed that when there is large imperfect items in minor(substitutable) products, the decision maker has to reduce its lot size whereas if there is a larger quantity of imperfect items in major(substituting) products, he has to increase its lot size. We have provided an effective algorithm to solve our model. Numerical example and sensitivity analysis are provided to test the validity and effectiveness of our model and to recommend some useful insights from a managerial viewpoint. The areas where our model can be applied include cement, metallic sheets, clothes, electronics items etc. This work can be extended to investigate a general multi-item inventory system. This model can be also modified to consider a problem with multi-level substitution in-spite of single level substitution. It should also be noticed that the system described in this paper assume continuous review inventory policy, which can be extended to a more complex policy in future work. One can incorporate the effect of inspection error, process shifting, learning and forgetting etc. in our model. In addition, we have assumed that the lead-times of the products are negligible in this study, which enabled us to work  with simplified model of the the aforementioned inventory system consisting of both products  and formulate the standard EOQ type equations. This condition can also be relaxed using positive or stochastic lead time in future. Another useful future modification can be to develop a model for a supply chain system with wholesaler and retailer with substitution  in a more general sense. More research is required with wide experimental data set to improve the present analysis herein and to obtain more generalized results, that can be applied to complex inventory systems.
\section*{Appendix-I}
We have obtained the expression for total average cost as follows:
\begin{eqnarray}\label{mm99}
\nonumber   TAC(\tau, T)&=&c_{h1}\Big[\frac{D_1+D_2}{2} T+{\frac { ( ( T-\tau ) D_2+T D_1 ) ^{2}E[p_1]}{ ( 1-E[p_1] ) ^{2}x_1T}}-{\frac {D_2{\tau}^{2}}{2 T}}
 \Big] +c_{h2} \Big[ {\frac {D_2{\tau}^{2}}{2 T}}+{\frac {E[p_2] D_2{\tau}^{2}}{ ( 1-E[p_2] ) ^{2}x_2T}} \Big] \\ && +{\frac {c_o}{T}}+D_2 c_t \Big[1-{\frac {\tau}{T}}\Big]
\end{eqnarray}
In order to prove the convexity of $TAC(\tau, T)$, we have to prove the following:\\
\vspace{0.1cm} \hspace{0.2cm}
$H(\tau,T)=\frac{\partial^2 TAC(\tau, T)}{\partial\tau^2}\frac{\partial^2 TAC(\tau, T)}{\partial T^2}-\frac{\partial^2 TAC(\tau, T)}{\partial\tau \partial T }>0$\\
\vspace{0.1cm} \hspace{0.2cm}
After calculation we obtain:\\
\vspace{0.1cm} \hspace{0.2cm}
$H(\tau,T)=\frac{1}{x_1 x_2 T^4(1-E[p_1])^2(1-E[p_2])^2}\Big[2 x_1 (1-E[p_1])^2\big[\big( -\frac{c_t^2}{2}-\tau(c_{h2}-c_{h1})(T-1)c_t D_2+c_o(c_{h2}-c_{h1})\big)(1-E[p_2])^2 x_2+2E[p_2]c_{h2}\big(\tau c_t(T-1)D_2-c_o\big)\big]+D_2x_2c_{h1}(1-E[p_2])^2(c_o-\tau c_t D_2(T-1))\Big]$
\\
\vspace{0.1cm} \hspace{0.2cm}
Consider an auxiliary function:\\
\vspace{0.1cm} \hspace{0.2cm}
$G(\tau,T)=2 x_1 (1-E[p_1])^2\big[\big( -\frac{c_t^2}{2}-\tau(c_{h2}-c_{h1})(T-1)c_t D_2+c_o(c_{h2}-c_{h1})\big)(1-E[p_2])^2 x_2 + 2 E[p_2]c_{h2}\big(\tau c_t(T-1)D_2-c_o\big)\big]+D_2x_2c_{h1}(1-E[p_2])^2(c_o-\tau c_t D_2(T-1))$\\
\vspace{0.2cm} \hspace{0.2cm}
We have to obtain the condition under which the function $G(\tau,T)>0$\\
\vspace{0.1cm} \hspace{0.2cm}
Now by inspection on the function G, we observe that if the condition $c_o>\frac{c_t^2}{2(c_{h2}-c_{h1})}+\tau(T-1)c_t D_2$ holds,
$G(\tau,T)$ will always positive. \\
\vspace{0.1cm} \hspace{0.2cm}
Then so will be $H(\tau,T)$. \\
Hence the proof of convexity of total cost follows.

\subsection*{}\label{refs}

\label{lastpage}

\end{document}